# An explicit substructuring method for overlapping domain decomposition based on stochastic calculus

Jorge Morón-Vidal, Francisco Bernal and Atsushi Suzuki

**Abstract**

In a recent paper [*F. Bernal, J. Morón-Vidal and J.A. Acebrón, Comp.& Math. App. 146:294-308 (2023)*] an hybrid supercomputing algorithm for elliptic equations has been put forward. The idea is that the interfacial nodal solutions solve a linear system, whose coefficients are expectations of functionals of stochastic differential equations confined within patches of about subdomain size. Compared to standard substructuring techniques such as the Schur complement method for the skeleton, the hybrid approach renders an explicit and sparse shrunken matrix—hence suitable for being substructured again. The ultimate goal is to push strong scalability beyond the state of the art, by leveraging the scope for parallelisation of stochastic calculus. Here, we present a major revamping of that framework, based on the insight of embedding the domain in a cover of overlapping circles (in two dimensions). This allows for efficient Fourier interpolation along the interfaces (now circumferences) and—crucially—for the evaluation of most of the interfacial system entries as the solution of small boundary value problems on a circle. This is both extremely efficient (as they can be solved in parallel and by the pseudospectral method) and free of Monte Carlo error. Stochastic numerics are only needed on the relatively few circles intersecting the domain boundary. In sum, the new formulation is significantly faster, simpler and more accurate, while retaining all of the advantageous properties of PDDSparse. Numerical experiments are included for the purpose of illustration.

## 1 Introduction

**Motivation.** The capability of efficiently solving large-scale boundary value problems (BVPs) motivated by realistic applications in science and technology is an increasingly important aspect of the digital and fourth industrial revolutions. In this paper we focus on linear elliptic BVPs—nonlinear ones are cast into an iteration thereof—and in two dimensions only (2D). The main reason is that the algorithm which will be presented is being developed from scratch, and the 2D case is easier to analyse, code and run, first. (Nonetheless, there are important "native" 2D applications, notably the Hele-Shaw as well as some atmospheric and oceanographic flows; see v.g. [4, 9] and references therein.)

Before going into our contribution, let us briefly describe the state of the art for massive linear systems stemming from discretisations of elliptic BVPs (see [12, 28, 30], among other). There are several options, but perhaps they can be exemplified by the Schur complement method [26]. The point is the "condensation" of all of the degrees of freedom (DoF) of the discretised algebraic system into the DoFs on the interfaces (or "skeleton") only, by taking the appropriate Schur complement of the original stiffness matrix. By leveraging the Schur complement, the skeletal linear system to be solved has orders of magnitude fewer DoFs than the original one (see [24] for an example). This is a huge simplification which comes at the price of sacrificing the original sparsity, and of not having an explicit matrix anymore. For the local inverses making up the Schur complement are almost never formed, owing to memory, accuracy, and cost considerations. The skeletal linear system is then solved on a parallel computer by a distributed and preconditioned Krylov iterative method, such as GMRES or BiCGStab [25]. The Schur complement structure involving local inverses is amenable to Krylov iterations, because the matrix-vector products underpinning the latter can be performed by solving local linear systems in parallel (on each iteration). On the other hand, keeping the number of iterations manageable calls for coarse-space preconditioners (see for instance [13] and references therein) which are themselves inherently sequential, and eventually set a cap to strong scalability [3] (i.e. the capability to efficiently use arbitrarily many independent computer processors working in parallel).

As the size and complexity of the original problems keeps growing due to the demand for realistic simulations, the Schur complement systems are becoming too big themselves, and ever more sophisticated preconditioners are required, in particular in the indefinite case [31]. Consequently, the number of cores in supercomputer facilities [1] is rapidly outpacing the software for using them (i.e. supercomputing algorithms), which improve rather incrementally.

**Previous work.** In the recent paper [7] we proposed PDDSparse. By inserting stochastic calculus (the Feynman-Kac formula specifically) into the formulation, it is able to produce an explicit, sparse, "one shot", substructured system in an embarrassingly parallel way. Concretely,

$$G\vec{u} = \vec{b},$$

where $\vec{u} = (u_1, \ldots, u_N)^\top$ are the nodal solutions on the $N$ interfacial *knots*[1] discretising the interfaces; and the components $(\vec{b})_k$ and nonzero entries $G_{ij}$ of the sparse matrix $G$ are certain expectations of SDEs, numerically integrated within small patches around each knot.[2] $G$ is about the same size as the Schur complement of the initial stiffness matrix, but the fact that $G$ is explicit and sparse paves the way for shrinking it again—for instance, by taking its Schur complement. The backward implication of this option of "shrinking twice" is that substantially bigger original systems might be tackled than it is possible now. However, the purpose of this article is another one—namely, constructing the skeletal system in a much faster and more accurate way than [7].

---

[1] We reserve the word "node" for supercomputer nodes.
[2] The knot's patch is the union of subdomains the knot belongs to.

**The new algorithm.** In this paper, we present a radical evolution of the PDDSparse algorithm. As will be discussed in detail in Section 2, the key modification consists in embedding the domain into a cover of *overlapping circles* (see Figure 1)—instead of nonoverlapping squares as in PDDSparse. Now, the circles are the subdomains, their circumferences are the interfaces, and the interfacial system is $C\vec{u} = \vec{r}$, but the nonzero entries of row $i$ corresponding to knot $\mathbf{x}_i$ *inside a circle* $\Omega_m$ *fully contained in the domain* are given by $C_{ij} = \mathcal{G}_j^m(\mathbf{x}_i)$ and $(\vec{r})_i = \mathcal{B}^m(\mathbf{x}_i)$, where the functions $\mathcal{G}_j^m(\mathbf{x})$ and $\mathcal{B}^m(\mathbf{x})$ obey well posed auxiliary BVPs on the circle. (In fact, all of the $\{\mathcal{G}_j^m\}$ functions for fixed $i$ stem from the same BVP with multiple right-hand sides, labeled by $j$.)

While the computation of the interfacial system remains as embarrassingly parallel as in PDDSparse, the advantage of the new method is that it is much faster and more accurate to solve BVPs on circles than it is to integrate SDEs inside them. This is due to the relatively poorer convergence rate of the error of numerical schemes for bounded diffusions [14, 21], as well as to the presence of statistical error.

The discussion above does not pertain to the perimeter subdomains, formed by the intersection of circles with the domain boundary, and whose corresponding entries of $C$ and $\vec{r}$ still must be estimated by Monte Carlo. But they are comparatively few.

Like PDDSparse's $G$, the new interfacial matrix $C$ has a rich structure [8]. A highly specific preconditioner is currently being designed which exploits it. Meanwhile—and in this paper—we are using a stopgap Schwarz type preconditioner. It is applied from a purely algebraic point of view, i.e, without reference to said properties of $C$ or to the underlying domain decomposition. Yet the fact that it can be adapted rather seamlessly attests to the potential for interconnection between our new algorithm, and the state of the art technology.

**Overview of the paper.** After this introduction, there are four more sections to this paper. In Section 2, the new algorithm is explained in detail. The stochastic numerics required for perimeter circles are the same as in [7], so we focus on the critical new part, i.e. Fourier interpolation along the circumferences, and the pseudospectral method for the local BVPs on the circles. Section 3 is devoted to various features of the new algorithm. Section 4 reports some initial numerical experiments. We also give details about our current (still preliminary) parallel implementation. Finally, Section 5 draws some conclusions and sets the guidelines for further work.

An Appendix with a self-contained exposition of the Feynman-Kac formulas for linear BVPs has been added in order to facilitate the reading to non-specialists, as that material tends to be scattered across the literature.

## 2 Formulation

We shall formulate the new algorithm for the two-dimensional elliptic BVP

$$\mathcal{L}u(\mathbf{x}) + c(\mathbf{x})u(\mathbf{x}) + f(\mathbf{x}) = 0 \text{ in } \Omega, \qquad u(\mathbf{x}) = g(\mathbf{x}) \text{ on } \partial\Omega, \tag{1}$$

where $\Omega \in \mathbb{R}^2$ is an open, bounded connected set; $c \leq 0$ in $\Omega$, and

$$\mathcal{L} = \frac{a_{xx}(x,y)}{2}\frac{\partial^2}{\partial x^2} + a_{xy}(x,y)\frac{\partial^2}{\partial x \partial y} + \frac{a_{yy}(x,y)}{2}\frac{\partial^2}{\partial y^2} + b_x(x,y)\frac{\partial}{\partial x} + b_y(x,y)\frac{\partial}{\partial y}, \tag{2}$$

such that the following matrix of second derivative coefficients is positive definite throughout $\Omega$—hence admitting a unique Cholesky decomposition:

$$[a_{ij}] := \begin{pmatrix} a_{xx} & a_{xy} \\ a_{xy} & a_{yy} \end{pmatrix} = \sigma\sigma^\top. \tag{3}$$

The reader is referred to the Appendix for details and references concerning the well posedness of (1), and to Section 3 for plausible extensions.

### 2.1 Discretisation and notation

The domain $\Omega$ is embedded into a cover of overlapping open balls $\{\bigcirc_m\}$

$$\bigcirc_m = \left\{\mathbf{x} \in \mathbb{R}^2 \,\middle|\, \|\mathbf{c}_m - \mathbf{x}\| < r_m\right\}, \qquad m = 1, \ldots, M \tag{4}$$

where $\mathbf{c}_m$ and $r_m$ are respectively the centre and the radius of the ball $\bigcirc_m$. In addition to fully covering $\Omega$, none of the balls can lie fully outside of it:

$$\Omega \subseteq \bigcup_{m=1}^{M} \bigcirc_m \quad \text{and} \quad \bigcirc_m \cap \Omega \neq \emptyset, \quad m = 1, \ldots, M. \tag{5}$$

In order to construct an overlapping domain decomposition $\Upsilon(\Omega) = \{\Omega_1, \Omega_2, \ldots, \Omega_M\}$, the subdomain $\Omega_k$ is defined as

$$\Omega_k = \bigcirc_k \cap \Omega. \tag{6}$$

Note that two kinds of subdomains arise: *floating* subdomains, which are circles fully immersed in $\Omega$, and *perimeter* subdomains, which intersect the boundary $\partial\Omega$ and may have an irregular, and typically nonsmooth, shape. Then

$$\Upsilon(\Omega) = \Upsilon_{\text{float}}(\Omega) \cup \Upsilon_{\text{per}}(\Omega), \tag{7}$$

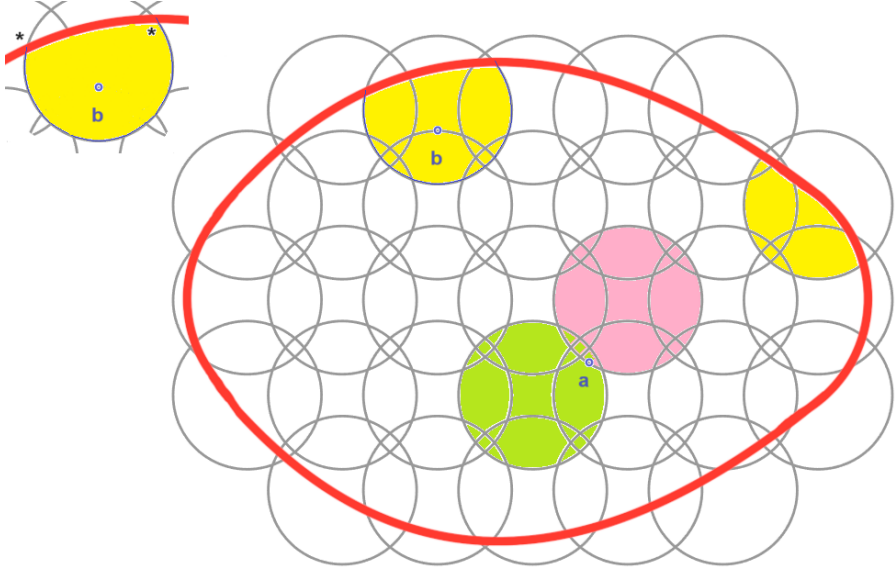

Figure 1: Discretisation into overlapping circles. The pink and green subdomains are floating, while the yellow ones are perimeter. Knot *a* belongs to the interface of the pink subdomain, and its $C_{ij}$ and $\vec{r}$ entries are computed by the pseudospectral method in the green subdomain. Those of knot *b* are obtained by stochastic simulation inside its perimeter subdomain. Outset: Note the kinks (*) at the junction between the boundary (red) and the perimeter arcs.

where $\Upsilon_{\text{float}}(\Omega)$ and $\Upsilon_{\text{per}}(\Omega)$ are the subsets of floating and perimeter subdomains, respectively. The *interface* $\Gamma_m$ of subdomain $\Omega_m$ is defined as

$$\Gamma_m = \partial \bigcirc_m \cap \overline{\Omega}, \qquad m = 1, \ldots, M. \tag{8}$$

Subdomain interfaces are circumferences of radius $r_m$, if $\Omega_m$ is floating, or circumference arcs, if $\Omega_m$ is perimeter. We make the assumption that $\partial \Omega$ is smooth enough that a cover of not too disparate circles can be constructed.

The interfaces are discretised into a set of $N$ points $\{\mathbf{x}_1, \ldots, \mathbf{x}_N\}$ called *knots*, which are the positions where the solution to the original BVP (1) will be computed. Due to the two different interpolation schemes (which will be explained in a moment), the knots must be equispaced on floating interfaces (full circumferences), but this is not strictly required on perimeter interfaces (arcs). Every knot must be in the interior of at least one subdomain. However, because of the circles overlap, a knot may (and often will) lie inside more than one subdomain. Finally, there must always be a knot at the junction between a perimeter arc and the actual boundary.

The index set $\Xi_m$ of the knots sitting on $\Gamma_m$ is called the interface stencil:

$$\Xi_m = \{i \mid \mathbf{x}_i \in \Gamma_m\}, \qquad m = 1, \ldots, M. \tag{9}$$

The knots in $\Xi_m$ are ordered in a particular way which will be explained in a moment. The stencil $\Xi_m$ and the interfaces $\Gamma_m$ are both called floating if $\Gamma_m$ is the interface of a floating subdomain. Otherwise, they are called perimeter.

Let #$S$ be the cardinal of a set $S$. Any smooth function $w : \Gamma_m \to \mathbb{R}$ can be linearly interpolated over the stencil of $\Gamma_m$ to arbitrary accuracy (within numerical precision) by refining the interface discretisation, i.e. by letting #$\Xi_m$ be large enough. More concretely,

$$\mathbf{z} \in \Gamma_m \quad \Rightarrow \quad w(\mathbf{z}) \approx \hat{w}(\mathbf{z}) \in span\left(\{w(\mathbf{x}_j)\}\,\big|\,j \in \Xi_m\right), \tag{10}$$

where $\hat{w}$ is the interpolator. We shall apply this ansatz to approximate the (still unknown) solution of the original BVP (1) along the interfaces, i.e. the $M$ interface-restricted functions $u|_{\Gamma_m} : \Gamma_m \to \mathbb{R}$. Note that improving the accuracy by refinement is possible because $\Gamma_m$ is a $C^\infty$ curve (a circumference or a circumference arc), and therefore $u|_{\Gamma_m}$ has the same smoothness as $u(\mathbf{x})$.

The natural way of parameterising a circumference arc is as

$$\Gamma_m = \left\{\mathbf{z} \,|\, \mathbf{z} = \mathbf{c}_m + r_m(\cos\theta, \sin\theta),\; \theta \in [\theta_m^{min}, \theta_m^{max})\right\}, \tag{11}$$

where $0 < \theta_m^{max} - \theta_m^{min} \leq 2\pi$. There is one single knot $\mathbf{x}_k$, with $1 \leq k \leq N$, associated to the $j^{th}$ knot on $\Gamma_m$, with $1 \leq m \leq M$ and $1 \leq j \leq$ #$\Xi_m$. Ordering the knots in $\Xi_m$ by increasing value of the local angle, we write

$$\theta_1 = \theta\!\left(\mathbf{x}_m^{(1)}\right) < \theta_2 = \theta\!\left(\mathbf{x}_m^{(2)}\right) < \ldots < \theta_{\#\Xi_m} = \theta\!\left(\mathbf{x}_m^{(\#\Xi_m)}\right) < \theta_m^{max}, \tag{12}$$

where the local index notation has been introduced as $\Gamma_m = \{\mathbf{x}_m^{(1)}, \ldots, \mathbf{x}_m^{(\#\Xi_m)}\}$.

Then, the sought-for interfacial nodal solution at position $\mathbf{z} \in \Gamma_m$ obeys

$$u(\mathbf{z}) \approx \hat{u}(\mathbf{z}) = \sum_{j=1}^{\#\Xi_m} H_j^m\!\left(\theta(\mathbf{z})\right) u\!\left(\mathbf{x}_m^{(j)}\right), \tag{13}$$

where the #$\Xi_m$ functions $H_j^m(\cdot)$ are the cardinal functions of the given one-dimensional interpolation scheme. Regardless of the latter, it is easy to check that cardinal functions have the Kronecker delta property: if $\theta_i, \theta_j$ are two nodes of the interpolation support, then

$$H_j^m(\theta_i) = \delta_{ij} \qquad \text{(1 if } i = j \text{ or 0 otherwise).} \tag{14}$$

As long as they are linear, we are in principle free to choose a different interpolation scheme for each interface. A natural choice is to use trigonometric interpolators for the full circumferences, and RBFs for the arcs.

## 2.2 Two families of interpolators and cardinal functions

**Interpolation on floating interfaces.** In this case, $\Gamma_m$ is a circumference. Any bounded function defined on $\Gamma_m$ will be periodic with period $2\pi$. Within this subsection, let us call $n = \#\Xi_m$. Since—by assumption—$u(\theta)$ is smooth enough and the knots are equally spaced, one can use the discrete Fourier transform to interpolate its value with high accuracy:

$$u|_{\Gamma_m} \approx \hat{u}(\theta) = \sum_{k=-n/2}^{n/2-1} \mathcal{A}_k \, e^{ik\theta}, \tag{15}$$

where $i = \sqrt{-1}$. We assume $n$ is even, so that $-n/2 \leq k \leq n/2 - 1$. (If $n$ was odd, then $-(n-1)/2 \leq k \leq (n-1)/2$.) The coefficients $\mathcal{A}_k$ are defined as

$$\mathcal{A}_k = \frac{1}{n} \sum_{j=1}^{n} u(\theta_j) \, e^{-ik\theta_j} \qquad k = \left\{-\frac{n}{2}, \ldots, \frac{n}{2}-1\right\}. \tag{16}$$

Then, (15) can be written as

$$\hat{u}(\theta) = \frac{1}{n} \sum_{k=-n/2}^{n/2-1} \sum_{j=1}^{n} u(\theta_j) \, e^{-ik\theta_j} \, e^{ik\theta}. \tag{17}$$

The imaginary part of (17) drops because $\hat{u}(\theta)$ is real, and

$$H_j^m(\theta) = \frac{1}{n} \sum_{k=-n/2}^{n/2-1} Re\left(e^{ik(\theta-\theta_j)}\right) = \frac{1}{n} \sum_{k=-n/2}^{n/2-1} \cos k(\theta - \theta_j). \tag{18}$$

**Remark 1** *Fourier cardinal functions (18) are periodic and $C^\infty$ over the complete interface of the floating subdomain. See Figure 2 (left).*

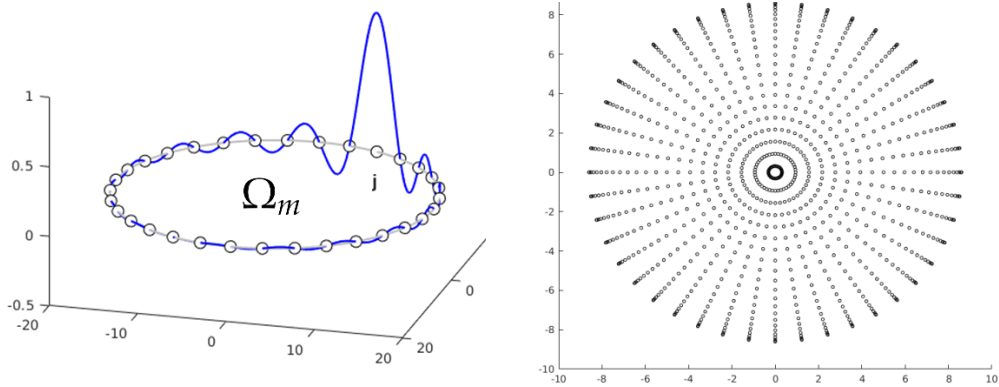

Figure 2: Discretisation of floating subdomains. (Left) The $H_j^m$ cardinal function (18). (Right) The Chebyshev grid (with $n_\theta$ and $n_r$ nodes along the angular and radial "directions", respectively) used by the pseudospectral method.

**Interpolation on perimeter interfaces.** The previous trigonometric interpolator is unsuitable when $\Gamma_m$ is a perimeter arc, since $u|_{\Gamma_m}$ is not a periodic function anymore. For perimeter interfaces, we have chosen radial basis function (RBF) interpolation, which attains spectral spatial accuracy for infinitely smooth functions on irregularly spaced stencils. The main drawback of RBF interpolation—namely ill-conditioned matrices with stencils of more that a few thousand knots—does not occur in this context. Specifically, we select the inverse multiquadric (IMQ) as RBF. The IMQ associated to knot $\mathbf{x}_m^{(j)} \in \Gamma_m$ is

$$\varphi_j(\theta) = \frac{1}{\sqrt{(\theta - \theta_j)^2 + c^2}}, \tag{19}$$

where $c^2$ is a tuneable parameter. Now we define the $n \times n$ symmetric matrix

$$K_{lk} = \varphi_k(\theta_l) \qquad 1 \leq l, k \leq n. \tag{20}$$

Following [7], the interpolator can then be expressed in Lagrange form as

$$\hat{u}(\theta) = \sum_{l=1}^{n} \sum_{k=1}^{n} u(\theta_k) \left(K^{-1}\right)_{lk} \varphi_l(\theta). \tag{21}$$

By comparison with (13), it is clear that the RBF cardinal functions are

$$H_j^m(\theta) = \sum_{l=1}^{n} \left(K^{-1}\right)_{lj} \varphi_l(\theta). \tag{22}$$

## 2.3 Interfacial linear system

The original BVP (1) restricted to the subdomain $\Omega_m$ reads

$$\begin{cases} \mathcal{L}u|_{\Omega_m}(\mathbf{x}) + c(\mathbf{x})u|_{\Omega_m}(\mathbf{x}) + f(\mathbf{x}) = 0 & \text{if } \mathbf{x} \in \Omega_m, \\ u|_{\Omega_m}(\mathbf{x}) = u|_{\Gamma_m}(\mathbf{x}) & \text{if } \mathbf{x} \in \Gamma_m, \\ u|_{\Omega_m}(\mathbf{x}) = g(\mathbf{x}) & \text{if } \mathbf{x} \in \partial\Omega \cap \partial\Omega_m, \end{cases} \tag{23}$$

where $u|_{\Gamma_m}(\mathbf{x})$ is the unknown interfacial solution. By ellipticity, $u|_{\Omega_m}$ coincides with the solution $u$ of (1) inside $\Omega_m$. The subdomain size is chosen in such a way that, upon discretisation, subdomain-restricted BVPs like (23) fit the memory of (and can be solved by) individual processors in the parallel computer.

The Dirichlet BC on the interface can now be approximated by its corresponding interpolator over the interface knots:

$$\begin{cases} \mathcal{L}u(\mathbf{x}) + cu(\mathbf{x}) + f(\mathbf{x}) = 0 & \text{if } \mathbf{x} \in \Omega_m, \\ u(\mathbf{x}) \approx \sum_{\mathbf{x}_j \in \Gamma_m} H_j^m(\mathbf{x}) u(\mathbf{x}_j) & \text{if } \mathbf{x} \in \Gamma_m, \\ u(\mathbf{x}) = g(\mathbf{x}) & \text{if } \mathbf{x} \in \partial\Omega \cap \partial\Omega_m. \end{cases} \tag{24}$$

Let knot $\mathbf{x}_i$ be inside $\Omega_m$ (not on its boundary). By the Feynman-Kac formula (check the Appendix), the pointwise solution to (24) at $\mathbf{x}_i \in \Omega_m$ is the expectation

$$u_i = \mathbb{E}\left[\left(\mathbb{1}_{\{\mathbf{X}_{\tau_m(\mathbf{x}_i)} \in \partial\Omega\}}\, g\left(\mathbf{X}_{\tau_m(\mathbf{x}_i)}\right) + \mathbb{1}_{\{\mathbf{X}_{\tau_m(\mathbf{x}_i)} \in \Gamma_m\}} \sum_{\mathbf{x}_j \in \Gamma_m} H_j^m\left(\mathbf{X}_{\tau_m(\mathbf{x}_i)}\right) u_j\right) Y_{\tau_m(\mathbf{x}_i)} + Z_{\tau_m(\mathbf{x}_i)} \,\bigg|\, \mathbf{X}_0 = \mathbf{x}_i\right]. \tag{25}$$

Note that we have written $u_i$ and $\{u_j\}$ rather than $u(\mathbf{x}_i)$ and $\{u(\mathbf{x}_j)\}$, for the solution of (23) and (24) is not exactly the same—although the latter converges to the exact one as the interfacial discretisation is refined and interpolation errors tend to zero.

In (25), the triple $(\mathbf{X}_{\tau_m(\mathbf{x}_i)}, Y_{\tau_m(\mathbf{x}_i)}, Z_{\tau_m(\mathbf{x}_i)})$ are obtained by integration of the SDE system (53) between $t = 0$ and $t = \tau_m(\mathbf{x}_i)$ (where $\sigma$ is given by (3)), and $\tau_m(\mathbf{x}_i)$ is the first exit time of $(\mathbf{X}_t)_{t \geq 0}$ from within $\Omega_m$, i.e.

$$\tau_m(\mathbf{x}_i) = \arg\min_{t \geq 0}\{\mathbf{X}_t \in \partial\Omega_m \,|\, \mathbf{X}_0 = \mathbf{x}_i\}, \tag{26}$$

which takes place at the first exit point $\mathbf{X}_{\tau_m(\mathbf{x}_i)} \in \partial\Omega_m$. Note that those statements are to made sense of on a trajectory basis, where a "trajectory" is a random realisation (continuous but nowhere differentiable) of the confined SDE for $(\mathbf{X}_t)_{0 < t \leq \tau_m(\mathbf{x}_i)}$. Reordering (25), we obtain

$$u_i - \sum_{\mathbf{x}_j \in \Gamma_m} \mathbb{E}\left[\mathbb{1}_{\{\mathbf{X}_{\tau_m(\mathbf{x}_i)} \in \Gamma_m\}} H_j^m\left(\mathbf{X}_{\tau_m(\mathbf{x}_i)}\right) Y_{\tau_m(\mathbf{x}_i)} \,\bigg|\, \mathbf{X}_0 = \mathbf{x}_i\right] u_j =$$
$$\mathbb{E}\left[\mathbb{1}_{\{\mathbf{X}_{\tau_m(\mathbf{x}_i)} \in \partial\Omega\}}\, g\left(\mathbf{X}_{\tau_m(\mathbf{x}_i)}\right) Y_{\tau_m(\mathbf{x}_i)} + Z_{\tau_m(\mathbf{x}_i)} \,\bigg|\, \mathbf{X}_0 = \mathbf{x}_i\right]. \tag{27}$$

Arranging the interfacial solution as the vector $\vec{u} = (u_1, \ldots, u_N)^\top$, and defining the $N \times N$ matrix $C$ and the $N \times 1$ vector $\vec{r}$, the interfacial nodal solutions solve the linear system

$$C\vec{u} = \vec{r} \tag{28}$$

whose coefficients are given by

$$C_{ij} = \begin{cases} 1 & \text{if } i = j, \\ -\mathbb{E}\left[\mathbb{1}_{\{\mathbf{X}_{\tau_m(\mathbf{x}_i)} \in \Gamma_m\}} H_j^m\left(\mathbf{X}_{\tau_m(\mathbf{x}_i)}\right) Y_{\tau_m(\mathbf{x}_i)} \,\bigg|\, \mathbf{X}_0 = \mathbf{x}_i\right] & \text{if } j \in \Xi_m, \\ 0 & \text{otherwise,} \end{cases} \tag{29}$$

and

$$(\vec{r})_i = \mathbb{E}\left[\mathbb{1}_{\{\mathbf{X}_{\tau_m(\mathbf{x}_i)} \in \partial\Omega\}}\, g\left(\mathbf{X}_{\tau_m(\mathbf{x}_i)}\right) Y_{\tau_m(\mathbf{x}_i)} + Z_{\tau_m(\mathbf{x}_i)} \,\bigg|\, \mathbf{X}_0 = \mathbf{x}_i\right]. \tag{30}$$

The nonzero columns of row $i$ of matrix $C$ are the indices $j$ of knots $\mathbf{x}_j$ sitting along the interfaces of the subdomain $\Omega_m$ that knot $\mathbf{x}_i$ belongs to the interior of. (The case when there are several of them is discussed later.) Therefore, $C$ is highly sparse when the number of subdomains $M$ is much larger that the typical number of knots along an interface. On the other hand, $C$ is nonsymmetric.

## 2.4 Computation of linear system coefficients

The expectations (29) and (30) ensure that the coefficients of the interfacial system can be computed by Monte Carlo stochastic simulation, and hence embarrassingly in parallel. However, thanks to the fact that most subdomains are circles, most coefficients can be computed orders of magnitude more efficiently.

**Coefficients of knots belonging to a floating subdomain.** In this case, the subdomain is not exposed to the actual boundary $\partial\Omega$, and all the trajectories from $\mathbf{x}_i \in \Omega_m$ end up on its interface. Since $\mathbf{X}_{\tau_m(\mathbf{x}_i)} \in \Gamma_m$ for all exit points, the corresponding nonzero off-diagonal elements of matrix $C$ become

$$C_{ij} = -\mathbb{E}\left[ H_j^m\left(\mathbf{X}_{\tau_m(\mathbf{x}_i)}\right) Y_{\tau_m(\mathbf{x}_i)} \middle| \mathbf{X}_0 = \mathbf{x}_i \right] \quad \text{if } \mathbf{x}_i \in \Omega_m \in \Upsilon_{\text{float}}. \tag{31}$$

Consider the auxiliary, subdomain-restricted (or *local*) family of BVPs

$$\begin{cases} \mathcal{L}\mathcal{G}_j^m(\mathbf{x}) + c(\mathbf{x})\mathcal{G}_j^m(\mathbf{x}) = 0 & \text{if } \mathbf{x} \in \Omega_m, \\ \mathcal{G}_j^m(\mathbf{x}) = -H_j^m(\mathbf{x}) & \text{if } \mathbf{x} \in \Gamma_m, \end{cases} \tag{32}$$

which differ from one another in the Dirichlet BC function, $-H_j^m(\mathbf{x})$. Since such functions are infinitely smooth (see Remark 1), (32) is well posed and admits a classical solution everywhere. In fact,

$$C_{ij} = \mathcal{G}_j^m(\mathbf{x}_i) \quad \text{if } \mathbf{x}_i \in \Omega_m \in \Upsilon_{\text{float}}. \tag{33}$$

Analogously, the expected value for the right-hand side (RHS) element associated with floating subdomains is

$$(\vec{r})_i = \mathbb{E}\left[ Z_{\tau_m(\mathbf{x}_i)} \middle| \mathbf{X}_0 = \mathbf{x}_i \right] \text{ if } \mathbf{x}_i \in \Omega_m \in \Upsilon_{\text{float}}, \tag{34}$$

which is the same as evaluating $\mathcal{B}^m(\mathbf{x}_i)$, where $\mathcal{B}^m$ solves the local BVP

$$\begin{cases} \mathcal{L}\mathcal{B}^m(\mathbf{x}) + c(\mathbf{x})\mathcal{B}^m(\mathbf{x}) + f(\mathbf{x}) = 0 & \text{if } \mathbf{x} \in \Omega_m, \\ \mathcal{B}^m(\mathbf{x}) = 0 & \text{if } \mathbf{x} \in \Gamma_m. \end{cases} \tag{35}$$

**Coefficients of knots belonging to perimeter subdomains.** In this case, the subdomain boundary contains both an interface (a circumference arc) and a portion of the actual boundary $\partial\Omega$. The BVP equivalent to (32) would be

$$\begin{cases} \mathcal{L}w_j(\mathbf{x}) + c(\mathbf{x})w_j(\mathbf{x}) = 0 & \text{if } \mathbf{x} \in \Omega_m, \\ w_j(\mathbf{x}) = -H_j^m(\mathbf{x}) & \text{if } \mathbf{x} \in \Gamma_m, \\ w_j(\mathbf{x}) = 0 & \text{if } \mathbf{x} \in \partial\Omega_m \cap \partial\Omega. \end{cases} \tag{36}$$

When discussing the discretisation, we specified that one knot must be placed at both intersections between the arc $\Gamma_m$ and the actual boundary $\partial\Omega$ (the asterisks in the outset of Figure 1). Let $\mathbf{x}_0$ be such a knot. By the Kronecker delta property of cardinal functions (14), $H_j^m(\mathbf{x}_0) = 0$ for all values of $j$, except for the specific cardinal function associated to knot $\mathbf{x}_0$, say $H_0^m$, for which $H_0^m(\mathbf{x}_0) = 1$. But on the actual boundary "side", the Dirichlet BC is 0. Consequently, the Dirichlet BC (and hence the solution $w_0$) undergoes a finite jump at $\mathbf{x}_0$—meaning that there is no classical solution $w_0(\mathbf{x})$ to (36) for that value of $j$. In other words, for perimeter subdomains, the expectation $\mathbb{E}\left[H_0^m\left(\mathbf{X}_{\tau_m(\mathbf{x}_i)}\right)Y_{\tau_m(\mathbf{x}_i)}\middle|\mathbf{X}_0 = \mathbf{x}_i\right]$ cannot be interpreted as the Feynman-Kac formula of a BVP.

Analogously, the perimeter equivalent to (35) would be

$$\begin{cases} \mathcal{L}w(\mathbf{x}) + c(\mathbf{x})w(\mathbf{x}) + f(\mathbf{x}) = 0 & \text{if } \mathbf{x} \in \Omega_m, \\ w(\mathbf{x}) = 0\,(\mathbf{x}) & \text{if } \mathbf{x} \in \Gamma_m, \\ w(\mathbf{x}) = g(\mathbf{x}) & \text{if } \mathbf{x} \in \partial\Omega_m \cap \partial\Omega. \end{cases} \quad (37)$$

whose Dirichlet BC is again discontinuous at the junctions $\Gamma_m \cap \partial\Omega$. Those matrix entries associated to perimeter patches can (and will) be computed by numerically integrating the corresponding SDEs. Since an ensamble of SDEs must be simulated in any case, it can produce all the nonzero $C_{ij}$ at once.

## 3 Discussion and implementation

### 3.1 Invertibility, stability and preconditioning of matrix $C$

In connection with PDDSparse, these topics are explored in depth in [8]. The main result is that the PDDSparse nodal matrix $G$ can be regarded as a perturbation of an asymptotic nonsingular M-matrix. The rich structure of the latter more than makes up for the lack of symmetry of $G$.

Notwithstanding the very different domain decomposition scheme, all of the proofs underpinning the theoretical analysis of PDDSparse and $G$ can be easily transplanted to the new algorithm and $C$—with one exception, most likely technical. While the cardinal functions of Fourier interpolation (18) are a nascent Dirac delta upon refinement (see Figure 2, left) as required by Theorem 4 in [8], the proof no longer works. The reason is that the latter relies on the nonoverlapping subdomain tessellation of PDDSparse, which does not carry over to the new overlapping decomposition in this paper. The claim of that proof—that the reduced matrix after deleting the knots on the boundary is irreducible—seems to hold empirically with $C$, though (check Table 1).

Therefore—pending a valid proof—it is likely that the new algorithm will also enjoy the advantageous M-matrixness-related properties of PDDSparse.

### 3.2 Choosing the floating subdomain for integration

In deriving (29) and (30), we implicitly assumed that there is one single subdomain that contains every knot. But owing to the circles overlap, a knot $\mathbf{x}_i$ will often belong to the interior of several subdomains. In that case, one can in principle choose any of them, or a weighted average. Our advice is to choose the subdomain whose interface is the farthest away from $\mathbf{x}_i$.[3]

Another possibility is calculating the Feynman-Kac expectations in all of those subdomains separately and taking the average. This has the drawback that increases the bandwidth of $C$ at row $i$, for it involves more stencil knots.

Nonetheless, for the sake of completeness we shall write next the formulas in this case. Letting $\Pi_i = \{\Omega_k \mid \mathbf{x}_i \in \Omega_k\}$,

$$C_{ij} = \begin{cases} \#\Pi_i & \text{if } i = j, \\ -\sum_{\Omega_m \in \Pi_i} \mathbb{E}\left[\mathbb{1}_{\{\mathbf{X}_{\tau_m(\mathbf{x}_i)} \in \Gamma_m\}} H_j^m\left(\mathbf{X}_{\tau_m(\mathbf{x}_i)}\right) Y_{\tau_m(\mathbf{x}_i)} \Big| \mathbf{X}_0 = \mathbf{x}_i\right] & \text{if } j \in \Gamma_m \\ 0 & \text{otherwise,} \end{cases} \quad (38)$$

and

$$(\vec{r})_i = \sum_{\Omega_m \in \Pi_i} \mathbb{E}\left[\mathbb{1}_{\{\mathbf{X}_{\tau_m(\mathbf{x}_i)} \in \partial\Omega\}} g\left(\mathbf{X}_{\tau_m(\mathbf{x}_i)}\right) Y_{\tau_m(\mathbf{x}_i)} + Z_{\tau_m(\mathbf{x}_i)} \Big| \mathbf{X}_0 = \mathbf{x}_i\right]. \quad (39)$$

### 3.3 The pseudospectral solver for floating subdomains

In order to take advantage of the floating subdomains being circles, and the solution assumed to be smooth, we resort to the pseudospectral collocation method to solve the local equations (32) and (35) on the floating subdomains. After changing to polar coordinates (see Figure 2, right) the method has been implemented based on [29, p.29], but using C++ instead of MATLAB.

After constructing the collocation matrix $A_k$ of a local BVP in $\Omega_k$, of order $n_\theta \times n_r$, the discrete solution $w$ on the pseudospectral nodes obeys

$$A_k w = (RHS), \quad (40)$$

which is solved with a direct solver—since local BVPs are meant to be small enough for single processors to tackle.

The nonzero matrix elements $C_{ij}$ on row $i$ (such that knot $\mathbf{x}_i \in \Omega_k$), as well as the skeletal RHS component $(\vec{r}_i)$, are the solution (evaluated on $\mathbf{x}_i$) to the $1 + \#\Xi_k$ BVPs (32) and (35), where the differential operators are the same—and thus the matrix $A_k$. This is equivalent to solving (32) and (35), with $1 + \#\Xi_k$ different RHSs. Thus, we perform upfront the QR factorisation of $A_k$ for each floating subdomain (namely $A_k = Q_k R_k$, where $R_k$ is an upper triangular matrix, and $Q_k$ is orthonormal), and then efficiently solve each of the BVPs as

$$R_k w_j = Q_k^\top (RHS)_j. \quad (41)$$

Notice, also, that the local BVPs are the same for every knot $\mathbf{x}_i$ inside the floating subdomain $\Omega_k$. Therefore, one single QR factorisation and $\#\Xi_k$ easy solves like (41) produce as many rows of $C$ as there are knots assigned to $\Omega_k$.

---

[3]The reason is that doing so is slightly more accurate and fosters the—yet to be proven—asymptotic M-matrixness of $C$, based on both experimental and theoretical observations.

## 3.4 Numerics for perimeter subdomains

Solving the system of SDEs (53) inside a perimeter subdomain takes the same kind stochastic numerics discussed in [7]. Let us stress the fact that the coefficients on perimeter subdomains can be calculated embarrassingly in parallel and taking full advantage of Graphics Processing Units (GPUs). (GPU-accelerated cores are an increasingly popular supercomputer architecture.)

## 3.5 Extensions

Elliptic BVPs with mixed BCs (of Dirichlet type on a portion of the boundary and of Robin/Neumann type on its complement) can be easily addressed, since the stochastic representation of such BVPs is well known (see the Appendix).

Mildly nonlinear BVPs can be tackled by leveraging the stochastic representation of the Fréchet derivative of the nonlinear operator [5, 23]. Three-dimensional domains might be discretised into overlapping balls or overlapping parallel cylinders.

**Remark 2** *Once the nodal system $C\vec{u} = \vec{r}$ has been built, it must be solved on a parallel computer. An optimised distributed solver which exploits the structure of C (with a tailored preconditioner) is a critical ingredient of the new algorithm, currently under development.*

# 4 Numerical experiments

The flow chart of the algorithm is on Figure 3. The code is written in C++ and parallelised with MPI and OpenMP. It was troubleshot and validated on Fugaku[4] [27], a supercomputer whose computing nodes are connected with a dedicated three-dimensional torus network, and each node consists of one multicore CPU A64FX based on ARM architecture. The CPU has 48 cores with a NUMA architecture and each group with 12 cores shares the memory channel. To perform hybrid parallel computation, a subset of cores in the NUMA unit will be assigned as shared memory computation by OpenMP. The master thread of OpenMP parallelisation is in charge of data communication with MPI.

**Stopgap linear solver for $C\vec{u} = \vec{r}$.** The distributed-memory solver of choice for an unsymmetric, general linear system such as $C\vec{u} = \vec{r}$ is GMRES [25]. For the time being, it is preconditioned with an off the shelf (algebraic) one-level restricted additive Schwarz (RAS) preconditioner [10], based solely on the adjacency graph associated to *C*. After constructing it, we perform an overlapping partition of the set of interfacial knots $\{\mathbf{x}_i\}_1^N$ into various subsets $\{w_p\}_1^P$, containing #$\omega_p$ elements each. The restriction of *C* to $w_p$ is defined as

$$C_p = R_p C R_p^\top, \tag{42}$$

---

[4] Our access grant allowed us to spend up to 48000 core-hours and 18432 cores simultaneously.

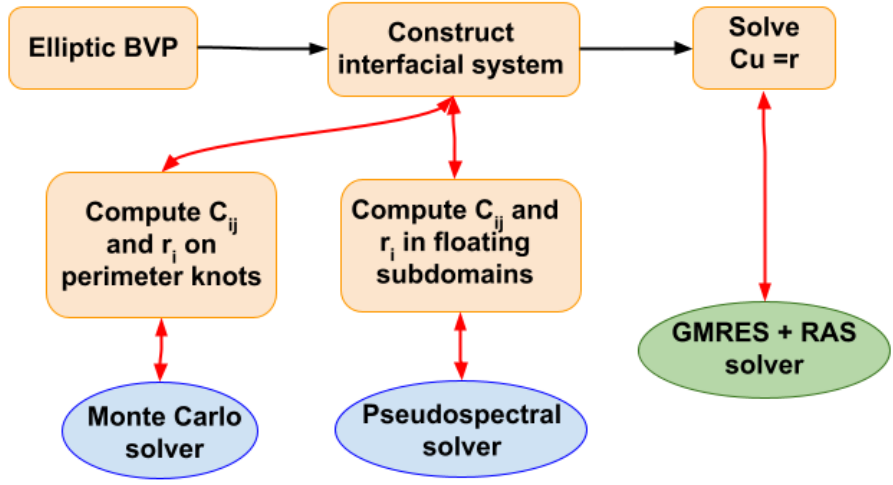

Figure 3: Squares contain the tasks and ovals, the tools used to accomplish them. Black/red arrows signify sequential/parallel execution.

where the $\#\omega_p \times N$ matrix $R_p$ is a restriction operator. The diagonal matrix $D_p$ is the discrete representation of the partition of unity restricted to each $w_p$,

$$(D_p)_{ii} = \frac{1}{\#\{q \mid \mathbf{x}_i \in \omega_q\}}, \qquad \sum_{p=1}^{P} R_p^\top D_p R_p = I_{N \times N}. \tag{43}$$

($I_{N\times N}$ is the identity matrix of order $N$.) Then, the RAS preconditioner is

$$M_{RAS}^{-1} = \sum_{p=1}^{P} R_p^\top D_p C_p^{-1} R_p. \tag{44}$$

GMRES solves the linear system with right preconditioner

$$CM_{RAS}^{-1}\left(M_{RAS}\vec{u}\right) = \vec{r} \tag{45}$$

as follows. Given $C$, we first perform an nonoverlapping decomposition of the graph of $C$ (using METIS [18]), after which adding neighbouring entries following the connectivity leads to $\{w_p\}_{p=1}^P$. Then, each core takes care of a subset of $C_p$ and performs the LU decomposition of those matrices (using MUMPS [2]). Finally, at each iteration of GMRES until convergence, those LU matrices are used to compute the action of $\{C_p^{-1}\}$ in (44) by forward/backward substitution. For simplicity, we have not yet enriched (44) with a coarse grid correction [19].

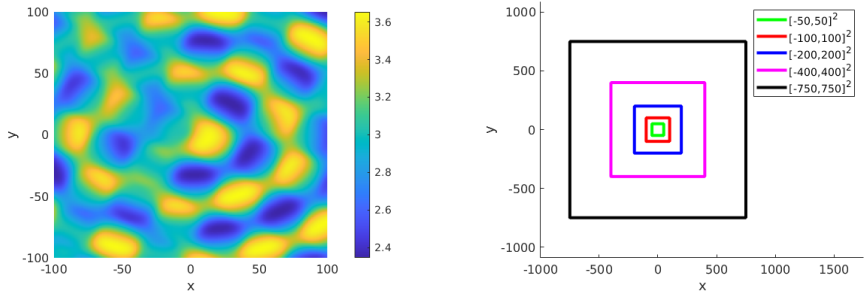

Figure 4: (Left) Exact solution of (46) in $[-100, 100]^2$. (Right) Test domains.

**The test problem.** We solve Poisson's equation $\nabla^2 u = F$ with exact solution

$$u_{ex}(x, y) = 3 + \tfrac{1}{3} \sin\left(\sqrt{1 + \tfrac{x^2}{100} + \tfrac{y^2}{50}}\right) + \tfrac{1}{3} \tanh\left[\sin\left(\tfrac{3x}{25} + \tfrac{y}{20}\right) + \sin\left(\tfrac{x}{20} - \tfrac{3y}{25}\right)\right]. \tag{46}$$

(Figure 4) and Dirichlet BCs equal to (46), on several squares: $[-50, 50]^2$, $[-100, 100]^2$, $[-200, 200]^2$, $[-400, 400]^2$, and $[-750, 750]^2$. In all the cases, the algorithm parameters are similar: the circles have radius between 8.55 and 9.43 and are arranged into a grid like that depicted on Figure 5 (left); there are 44 knots per circumference, arranged like in Figure 2 (left). On floating subdomains, there are also $n_\theta = 44$ nodes in the angular "direction" (coincident with the interfacial knots), and $n_r = 22$ along the radius for a total of 968 pseudospectral nodes per subdomain (Figure 2, right). As a reference, the pseudospectral collocation matrix has a condition number of order $10^5$ and computes the solution of (46) inside a subdomain with a maximum error of order $10^{-12}$. On perimeter subdomains, we used the Gobet-Menozzi method [15], pathwise variance reduction as in [6], $\tilde{N} = 5000$ trajectories per perimeter knot, and timestep $h = 0.015$. The shape parameter in (19) was $c^2 = 1.5$. We stress that those parameters were chosen rather casually.

Table 1 lists an output selection resulting from these simulations. The fact that both the spectral radius of $C - I$ is almost below one, and $\max_{1 \le i \ne j \le N} C_{ij}$ is almost nonpositive suggests that $C$ may be regarded as a perturbed M-matrix [17].

| domain | $[-50, 50]^2$ | $[-100, 100]^2$ | $[-200, 200]^2$ | $[-400, 400]^2$ | $[-750, 750]^2$ |
|---|---|---|---|---|---|
| subdomains | 100 | 400 | 1600 | 6400 | 22500 |
| N | 3636 | 15997 | 67157 | 275077 | 978036 |
| sparsity | 2.066% | 0.502% | 0.123% | 0.030% | 0.009% |
| $\rho(C - I)$ | 0.9558 | 0.9900 | 0.9976 | 0.9994 | 0.9998 |
| $\max_{1 \leq i \neq j \leq N} C_{ij}$ | 0.0445 | 0.0767 | 0.0800 | 0.0801 | 0.0494 |
| $\kappa_2(C)$ | 84 | 371 | 1560 | 6390 | 22700 |
| RMS error | $3.72 \times 10^{-5}$ | $3.01 \times 10^{-5}$ | $1.86 \times 10^{-5}$ | $1.56 \times 10^{-5}$ | $1.19 \times 10^{-5}$ |
| GMRES iterations | 19 | 33 | 51 | 73 | 101 |
| solve $C\vec{u} = \vec{r}$ (s.) | 1.55 | 4.2 | 17.3 | 85.2 | 370 |
| total time (s.) | 76 | 203 | 655 | 2171 | 7356 |

Table 1: Using 24 nodes, 96 MPI processes and 1152 cores (no GPUs).

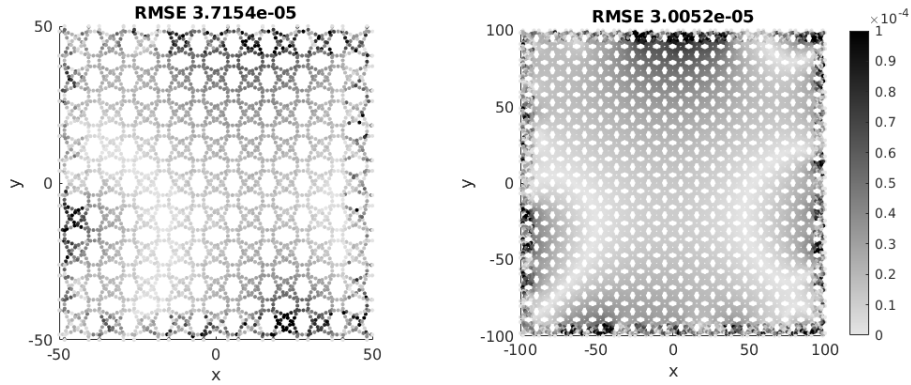

Figure 5: Error on the interfacial knots in domains $[-50, 50]^2$ and $[-100, 100]^2$.

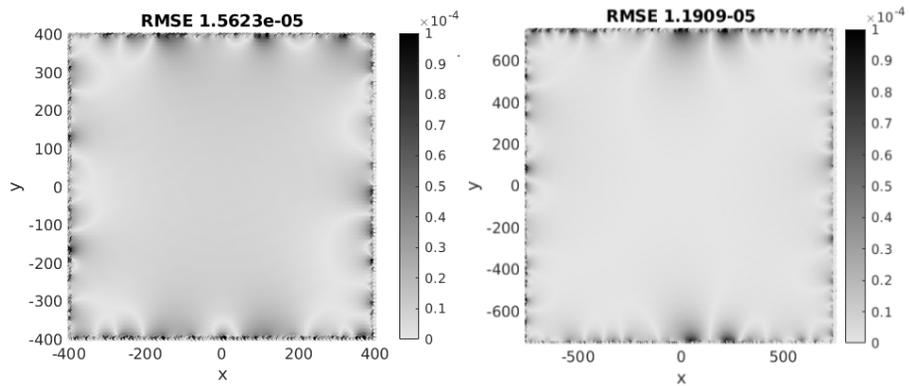

Figure 6: Nodal error in $[-400, 400]^2$ and $[-750, 750]^2$ (beware different scales).

Figures 5 and 6 show the nodal error throughout the domain(s). It can be observed that:

i) The error is larger close to the boundary. This is a consequence of the relatively low accuracy and statistical error of Monte Carlo simulations.

ii) Well inside the square domains, the error is much lower, fairly homogeneous, and independent of the size of the square.

iii) Still, the error in the interior is larger than what would be expected from pseudospectral simulations alone, indicating that some of the error is propagated inwards from the boundary by the nodal system $C\vec{u} = \vec{r}$.

iv) The new algorithm is $10^2 - 10^3$ times more accurate than PDDSparse.

While the accuracy of the algorithm is quite satisfactory, the admittedly complex structure of the error makes it difficult to disentangle the effect of its various sources. Nonetheless, the aggregate error does converge roughly as expected with $h$ and $\tilde{N}$. For instance, with $N = 978036$ and $M = 150^2$, the triples $(h, \tilde{N}, \text{RMS error})$ are: $(0.015, 5000, 1.2\times 10^{-5})$; $(0.005, 2500, 2.0\times 10^{-6})$; and $(0.00167, 2500, 1.0 \times 10^{-6})$. This is reasonable given the perimeter subdomain kinks (Figure 1), and the discontinuous payoff.

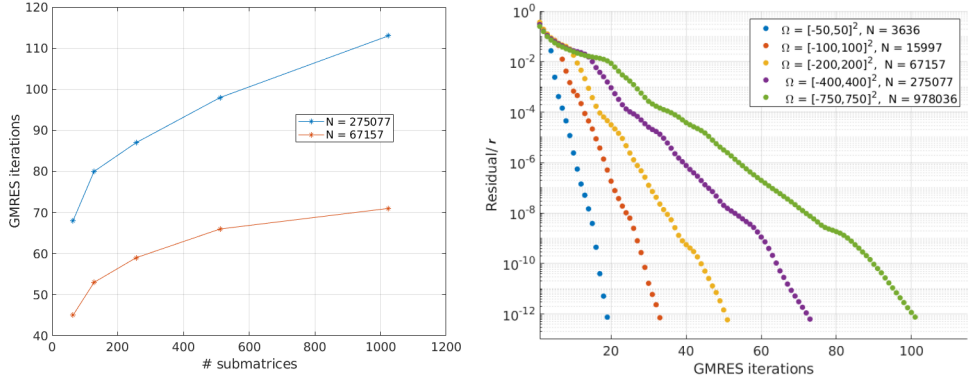

Figure 7: The curves on the left are fitted to $y = a + x^\alpha$ with $\alpha = 0.18, a = 3.5$ (blue) and $\alpha = 0.16, a = 3.2$ (orange). Right: $r = \|\vec{r}\|_2$.

As expected for single-level RAS, the number of normalised[5] GMRES iterations for $C\vec{u} = \vec{r}$ is not scalable—see Figure 7 (left). However, it grows with a lower power $\alpha$ of the number of METIS submatrices than the theoretical rate in standard domain decomposition: about $\alpha = 0.2$ instead of $\alpha = 0.5$ [30]. (This is probably due to the significantly larger submatrix coupling of our method.) Let us also point out that the number of GMRES iterations grows sublinearly with the discretisation size (Figure 7, right), and much more slowly than the condition number (see Table 1).

---

[5]Defined as the $\ell_2$ norm of the residual over that of $\vec{r}$.

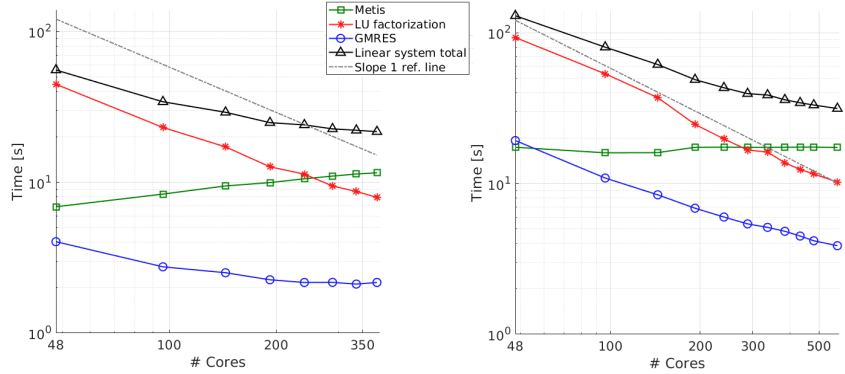

Figure 8: Empirical strong scalabitity of the linear solver, by component.

For completeness, the time taken by the various tasks of (45) is dissected in Figure 8. On the right, $P = 960$ is fixed. It scales better (in the strong sense), but is overall slower than aligning $P$ to the number of MPI processes[6] (left). The bottleneck is METIS. Without taking it into account, strong scaling of this stopgap implementation would be improved by larger $N$, and a coarse space.

# 5 Conclusions and future work

The algorithm presented here offers a fast, accurate and embarrassingly parallel way of constructing a skeletal matrix of the original problem which is both sparse and explicit, in contrast to the state of the art substructuring methods. The latter can now be applied to the reduced order system generated by the new algorithm, which in turn spells significantly larger original problems. A tailored Schur complement condensation of $C\vec{u} = \vec{r}$ is, in fact, our next objective.

# Acknowledgements

FB and JMV were funded by grants 2018-T1/TIC-10914 and 2022-5A/TIC-24233 of the Madrid Regional Government in Spain, RIKEN R-CCS Fugaku Trial Access call HP220345, and partially by Spanish AEI grant PID2020-115088RB-I00. JMV thanks the hospitality of RIKEN Center for Computational Science in Kobe (Japan).

---

[6] In those simulations, an MPI process executes an OpenMP computation with 3 cores.

# A  Stochastic representation of linear second order BVPs

Let $d \geq 2$, $\Omega \subset \mathbb{R}^d$ be a bounded domain, and $\Omega = \overline{\Omega} \cup \partial\Omega$, where the open set $\overline{\Omega}$ is the interior of the domain and $\partial\Omega$ its boundary. Consider the linear parabolic BVP of second order with mixed BCs:

$$\begin{cases} \frac{\partial u}{\partial t} = \mathcal{L}(\mathbf{x},t)u + c(\mathbf{x},t)u + f(\mathbf{x},t) & \text{if } 0 < t \leq T, \mathbf{x} \in \overline{\Omega}, \\ u = p(\mathbf{x}) & \text{if } t = 0, \mathbf{x} \in \Omega, \\ u = g(\mathbf{x},t) & \text{if } 0 < t \leq T, \mathbf{x} \in \partial\Omega_A, \\ \frac{\partial u}{\partial N} = \varphi(\mathbf{x},t)u + \psi(\mathbf{x},t) & \text{if } 0 < t \leq T, \mathbf{x} \in \partial\Omega_R, \end{cases} \quad (47)$$

where $T > 0$, $\varphi(\mathbf{x},t) \leq 0$, and

$$\mathcal{L}(\mathbf{x},t)u = \frac{1}{2}\sum_{i,j=1}^{d} a_{ij}(\mathbf{x},t)\frac{\partial^2 u}{\partial x_i \partial x_j} + \sum_{k=1}^{d} b_k(\mathbf{x},t)\frac{\partial u}{\partial x_k}. \quad (48)$$

($\mathcal{L}$ is called the differential generator in the SDE context.) The matrix $A(\mathbf{x},t) := [a_{ij}]_{i,j=1}^{d}$ is positive definite, and $\mathbf{b}(\mathbf{x},t) := (b_1,\ldots,b_d)^\top$ is called the drift. All of the coefficient functions in (47), namely $a_{ij}, b_i, c, f, p, g, \varphi$ and $\psi$ are assumed continuous, and complying with the compatibility conditions at time $t = 0$ (see [20] or [11, formulas (2.15)-(2.17)]). The boundary is decomposed as $\partial\Omega = \partial\Omega_A \cup \partial\Omega_R \cup \partial\Omega_S$, such that $\partial\Omega_A \cap \partial\Omega_R = \partial\Omega_A \cap \partial\Omega_S = \partial\Omega_R \cap \partial\Omega_S = \emptyset$. (Note that $\partial\Omega_A$ and $\partial\Omega_R$ are known as absorbing and reflecting boundaries in the SDE literature.) The *outward* unit normal vector $\mathbf{N}$ is assumed to be well defined on the boundary save perhaps on a set $\partial\Omega_S$; $\partial\Omega_A$ stands for the portion of the boundary (if any) where Dirichlet BCs are imposed; and on $\partial\Omega_R$, BCs involving the normal derivative, (i.e. $\mathbf{N}^\top \nabla u$) hold, where $\nabla u = (\partial u/\partial x_1, \ldots, \partial u/\partial x_d)^\top$. (BCs involving oblique, rather than normal, derivatives—the so-called *third boundary value problem*—will not be considered in this paper.) Such BCs are of Neumann type iff $\varphi = 0$, or of Robin type otherwise.

**Sufficient conditions for existence of a unique classical solution to (47).** If $\partial\Omega \in C^2$ (thus $\partial\Omega_S = \emptyset$) and $\partial\Omega_R = \emptyset$ (respectively $\partial\Omega_A = \emptyset$), theorem 2.6 (respectively theorem 2.7) in [11] (see also [20]) ensure the existence of a unique classical solution with regularity dependent on that of the BVP coefficients and of $\partial\Omega$. (By a classical solution, we mean $u(\mathbf{x},t)$ which lives in the Hölder space $C^{1,2}([0,T] \times \overline{\Omega})$.) When both $\partial\Omega_A \neq \emptyset \neq \partial\Omega_R$, this connection is less general and more dependent on the smoothness of the boundary [22].

**Further conditions for the stochastic representation.** Let $\sigma(\mathbf{x},t)$ (called the diffusion matrix) be defined by $\sigma(\mathbf{x},t)\sigma^\top(\mathbf{x},t) = A(\mathbf{x},t)$ (this is always possible since $A$ is positive definite). The following result is an extension of the well-known Feynman-Kac formulas for purely reflected ($\partial\Omega_A = \emptyset$) and purely absorbed ($\partial\Omega_R = \emptyset$) diffusions, adapted from [11, theorem 2.5] and [21, chapter 6].

**Theorem 1** *Assume that: i) a classical unique solution to (47) does exist; ii) there exists a constant $L > 0$ such that for $\mathbf{x}, \mathbf{y} \in \overline{\Omega}$ and $t \in [0, T]$*

$$\|\sigma(\mathbf{x}, t) - \sigma(\mathbf{y}, t)\| \leq L\|\mathbf{x} - \mathbf{y}\|, \tag{49}$$

$$\|\mathbf{b}(\mathbf{x}, t) - \mathbf{b}(\mathbf{y}, t)\| \leq L\|\mathbf{x} - \mathbf{y}\|; \tag{50}$$

*iii) $\partial\Omega$ is piecewise $C^1$ (i.e. $C^1$ save on maybe a set $\partial\Omega_S$); and iv) either $\overline{\Omega}$ is convex, or $u(\mathbf{x}, t)$ can be extended to a function $C^{1,2}([0, T] \times \mathbb{R}^d)$. Then, for $0 \leq t \leq T$ the pointwise solution of (47) admits the following stochastic representation:*

$$u(\mathbf{x}_0, t) = \mathbb{E}[\phi_\tau] := \mathbb{E}\left[q(\mathbf{X}_\tau)Y_\tau + Z_\tau\right], \tag{51}$$

*where*

$$q(\mathbf{X}_\tau) = \begin{cases} g(\mathbf{X}_\tau, T - \tau), & \text{if } \tau < T, \\ p(\mathbf{X}_T), & \text{if } \tau \geq T, \end{cases} \tag{52}$$

*and the processes $(\mathbf{X}_t, Y_t, Z_t, \xi_t)$ are governed by the following set of SDEs driven by a standard $d$-dimensional Wiener process $\mathbf{W}_t$:*

$$\begin{cases} d\mathbf{X}_t = \mathbf{b}(\mathbf{X}_t, T - t)dt + \sigma(\mathbf{X}_t, T - t)d\mathbf{W}_t - \mathbf{N}(\mathbf{X}_t)d\xi_t & \mathbf{X}_0 = \mathbf{x}_0, \\ dY_t = c(\mathbf{X}_t, T - t)Y_t dt + \varphi(\mathbf{X}_t, T - t)Y_t d\xi_t & Y_0 = 1, \\ dZ_t = f(\mathbf{X}_t, T - t)Y_t dt + \psi(\mathbf{X}_t, T - t)Y_t d\xi_t & Z_0 = 0, \\ d\xi_t = \mathbb{1}_{\{\mathbf{X}_t \in \partial\Omega_R\}}dt & \xi_0 = 0. \end{cases} \tag{53}$$

*Above, $\mathbb{1}_{\{H\}}$ is the indicator function (1 if H is true and 0 otherwise); $d(\cdot, \cdot)$ denotes the Euclidean distance from a point to a set; $\tau = \inf_t\{\mathbf{X}_t \in \partial\Omega_A\}$ is the "first exit time" (or first passage time) from $\Omega$; which takes place at the "first exit point" $\mathbf{X}_\tau \in \partial\Omega_A$; and $\xi_t$ is the "local time" up to time $t$.*

Intuitively, $\xi_t$ is the amount of time (up to time $t$) that $\mathbf{X}_t$ spends infinitesimally close to $\partial\Omega_R$; rigorously the pair of processes $(\mathbf{X}_t, \xi_t)$ are the solution of Skorohod's problem inside a bouded domain with absorbing-normally reflecting boundary—see [11, section 2] for details. Sometimes it is more convenient to use the alternative formula

$$\phi_\tau = \mathbb{1}_{\{\tau \geq T\}} p(\mathbf{X}_T)\Phi(T) + \mathbb{1}_{\{\tau < T\}} g(\mathbf{X}_\tau, T - \tau)\Phi(\tau) + \int_0^{\min(T,\tau)} f(\mathbf{X}_t, T - t)\Phi(t)dt$$

$$+ \int_0^{\min(T,\tau)} \psi(\mathbf{X}_t, T - t)\Phi(t)d\xi_t, \tag{54}$$

where

$$\Phi(t) = \exp\left(\int_0^t c(\mathbf{X}_s, s)ds + \int_0^t \varphi(\mathbf{X}_s, s)d\xi_s\right). \tag{55}$$

**Elliptic equations.** Equation (47) can be formally transformed into an elliptic BVP with mixed BCs by assuming that $\partial u/\partial t = 0$, thus dropping the dependence on time from $u$ and all the coefficients; letting $T \to \infty$; and dropping the initial condition $p$. If $\partial \Omega_R = \emptyset$ (purely stopped diffusions / Dirichlet BCs), the stochastic representation derived from Theorem 1 still holds as long as $c \leq 0$ and $\mathbb{E}[\tau] < \infty$ [16, chapter 4.4.5]. (Note that in the purely reflected case, the latter condition would be impossible.) To the best of our knowledge, there is no rigorously proven stochastic representation for elliptic BVPs with mixed BCs. However, we will assume in the remainder of this paper that, if: i) the time-independent equivalent conditions of those in Theorem 1 are in place; ii) $c \leq 0$; and iii) $\mathbb{E}[\tau] < \infty$; then the time-independent equivalent representation given by Theorem 1 holds. Note that in that case, $\mathbb{1}_{\{\tau \geq T\}} p(\mathbf{X}_T)\Phi(T) = 0$ in (54) and $0 < \Phi(t) \leq 1$ in (55). The appropriate stochastic representations are

$$u(\mathbf{x}) = \mathbb{E}\Big[g(\mathbf{X}_\tau)e^{\int_0^\tau c(\mathbf{X}_s)ds + \int_0^\tau \varphi(\mathbf{X}_s)d\xi_s} + \int_0^\tau f(\mathbf{X}_t)e^{\int_0^t c(\mathbf{X}_s)ds + \int_0^t \varphi(\mathbf{X}_s)d\xi_s}dt$$

$$+ \int_0^\tau \psi(\mathbf{X}_t)e^{\int_0^t c(\mathbf{X}_s)ds + \int_0^t \varphi(\mathbf{X}_s)d\xi_s}d\xi_t\Big] = \mathbb{E}\Big[g(\mathbf{X}_\tau)Y_\tau + Z_\tau\Big], \quad (56)$$

where $(\mathbf{X}_s, Y_s, Z_s, \xi_s)$ are the solution of the equivalent system to (53) resulting from dropping the dependence on $T - t$ in all the coefficients.